\documentclass[reqno,10pt]{amsart}
\usepackage{amsmath, amsfonts}

\usepackage{color}
\usepackage{amsmath}
\usepackage{amsfonts}
\usepackage{amssymb}
\usepackage{amsthm}
\usepackage{newlfont}
\usepackage{graphicx}
\usepackage{hyperref}

\newtheorem{thm}{Theorem}[section]
\newtheorem{cor}[thm]{Corollary}
\newtheorem{lem}[thm]{Lemma}
\newtheorem{prop}[thm]{Proposition}
\newtheorem{conj}[thm]{Conjecture}
\theoremstyle{definition}
\newtheorem{defn}[thm]{Definition}
\theoremstyle{remark}
\newtheorem{rem}[thm]{Remark}

\numberwithin{equation}{section}
\newcommand{\eps}{\varepsilon}

\newcommand{\lsm}{\lesssim}

\newcommand{\C}{{\mathbb{C}}}
\newcommand{\R}{{\mathbb{R}}}

\newcommand{\ltrd}{L^2_x (\R^d)}

\newcommand{\pn}{P_N}

\newcommand{\ins}{[0,\frac 1{\sqrt N}]}
\newcommand{\insd}{[0,\frac 1{\sqrt N}]\times \R^d}
\newcommand{\sqn}{\frac 1{\sqrt N}}

\newcounter{smalllist}
\newenvironment{SL}{\begin{list}{{$($\roman{smalllist}\/$)$\hss}}{%
\setlength{\topsep}{0mm}\setlength{\parsep}{0mm}\setlength{\itemsep}{0mm}%
\setlength{\labelwidth}{2.0em}\setlength{\itemindent}{2.5em}\setlength{\leftmargin}{0em}\usecounter{smalllist}%
}}{\end{list}}

\title[Regularity of APMS solutions]{Regularity of almost periodic modulo scaling solutions
for mass-critical NLS and applications}
\author{Dong Li}
\address{Department of Mathematics, University of Iowa, 14 MacLean Hall, Iowa City, IA, 52240}
\author{Xiaoyi Zhang}
\address{Academy of Mathematics and System Sciences, Beijing, and Department of Mathematics, University of
Iowa, 14 MacLean Hall, Iowa City, IA, 52240}


\begin{document}

\maketitle

\begin{abstract}

In this paper, we consider the $L_x^2$ solution $u$ to mass critical
NLS $iu_t+\Delta u=\pm |u|^{\frac 4d} u$. We prove that in
dimensions $d\ge 4$, if the solution is spherically symmetric and is
\emph{almost periodic modulo scaling}, then it must lie in
$H_x^{1+\eps}$ for some $\eps>0$. Moreover, the kinetic energy of
the solution is localized uniformly in time. One important
application of the theorem is a simplified proof of the scattering
conjecture for mass critical NLS without reducing to three enemies
\cite{ktv:2d}, \cite{kvz:blowup}. As another important application,
we establish a Liouville type result for $L_x^2$ initial data with
ground state mass. We prove that if a radial $L_x^2$ solution to
focusing mass critical problem has the ground state mass and does
not scatter in both time directions, then it must be global and
coincide with the solitary wave up to symmetries. Here the ground
state is the unique, positive, radial solution to elliptic equation
$\Delta Q-Q+Q^{1+\frac 4d}=0$. This is the first rigidity type
result in scale invariant space $L_x^2$.

\end{abstract}

%\tableofcontents
\section{Introduction}
\subsection{Main results}
We consider the $d$-dimensional mass critical nonlinear
Schr\"odinger equation
\begin{align}\label{nls}
iu_t+\Delta u=\mu |u|^{\frac 4d} u=:F(u).
\end{align}
Here, $\mu=\pm1$ with $\mu=+1$ known as the "defocusing" and
$\mu=-1$ as the "focusing" case. The name "mass-critical" refers to
the fact that the scaling symmetry
\begin{align}\label{scaling}
u(t,x)=\lambda^{\frac d2}u(\lambda^2 t,\lambda x),
\end{align}
leaves both the equation and the mass invariant. Here the mass is
defined as
\begin{align}\label{mass}
\mbox{Mass: } M(u(t))=\int_{\R^d} |u(t,x)|^2 dx=M(u_0).
\end{align}

The precise meaning of the solution we discuss throughout the paper
is the following
\begin{defn}[Solution] A function $u:I\times\R^d\to\C$ on a
non-empty time interval $I\subset \R$ is a strong $L_x^2(\R^d)$
solution (or solution for short) if it lies in the class
$C_t^0L_x^2(K\times \R^d)\cap L_{t,x}^{\frac{2(d+2)}d}(K\times\R^d)$ for all
compact $K\subset I$, and we have the Duhamel formula
\begin{align}
u(t_1)=e^{i(t_1-t_0)\Delta}u(t_0)-i\int_{t_0}^{t_1}
e^{i(t_1-t)\Delta}F(u(t))dt
\end{align}
for all $t_0,t_1\in I$. Here $e^{it\Delta}$ is the propagator for
free Schr\"odinger equation. We say that $u$ is a maximal-lifespan
solution if the solution can not be extended to any strictly larger
interval. We say that $u$ is global if $I=\R$.
\end{defn}

The standard local theory for the above defined solution was worked
out by Cazenave and Weissler \cite{caz}. They constructed the local
in time solution for arbitrary initial data in $L_x^2(\R^d)$. They
also showed that the solution depends continuously on the initial
data in the same space. However due to the criticality of the
problem, the lifespan of the local solution depends on the profile
of the initial data instead of the mere $L_x^2$-norm. When the initial
data is small enough, they proved the solution exists globally and
scatters in the following sense: there exist unique $u_{\pm}\in
L_x^2(\R^d)$ such that
\begin{align}\label{defn:scat}
\lim_{t\to\infty}\|u(t)-e^{it\Delta }u_{+}\|_{L_x^2}=\lim_{t\to -\infty}
\|u(t)-e^{it\Delta} u_{-}\|_{L_x^2}
=0.
\end{align}

Whilst the local theory is fairly complete, the understanding of the
global theory for large solutions is still only partial. Briefly
speaking, the global theory for large solutions amounts to proving
the global wellposeness and scattering for generic $L_x^2$ initial
data in the defocusing case; investigating the long time behavior of
global solutions, characterizing the structure and profile of finite
time blowup solutions in the focusing case and so on. In recent
years, by using concentration compactness tools developed and used
in
\cite{BegoutVargas,keme06s,keraani,keraani-2,klvz,ktv:2d,kvz:blowup,lz:2d,lz:six,merle_duke},
one can address part of these problems by exploring the properties
of a large class of solutions which have certain compactness
properties. To this end, following \cite{compact}, we introduce

\begin{defn}[Almost periodic modulo symmetry solutions]\label{amps}
Let $u$ be the maximal-lifespan solution of \eqref{nls} on time
interval $I$. Let $I_0\subset I$ be a subinterval. We say $u$ is
\emph{almost periodic modulo symmetries} on $I_0$ if there exists functions
$x(t), N(t), \xi(t), \theta(t)$ with $t\in I_0$ such that the orbit
\begin{align*}
\biggl\{e^{i\theta(t)}e^{ix\cdot\xi(t)}N(t)^{-\frac
d2}u\biggl(t,\frac{x-x(t)}{N(t)}\biggr),\ t\in I_0\biggr\}
\end{align*}
is precompact in $L_x^2(\R^d)$. By Arzela-Ascoli Theorem, an equivalent way to write this definition is the following:
there exists a function $C(\eta)$ such that for any $\eta>0$,
\begin{align*}
\int_{|x-x(t)|>\frac{C(\eta)}{N(t)}}|u(t,x)|^2 dx\le \eta,\\
\int_{|\xi-\xi(t)|>C(\eta) N(t)}|\hat u(t,\xi)|^2 d\xi \le \eta.
\end{align*}
In particular, we call $u$ is
\emph{almost periodic modulo scaling} on $I_0$ if in the above
definition, $x(t)=\xi(t)\equiv 0$ for all $t\in I_0$.
\end{defn}

In the above definition, the parameter $N(t)$ is the frequency
scale. In the physical space, its reciprocal corresponds to the
concentration size of the solution. The parameter $x(t)$, $\xi(t)$
correspond to the center of mass at physical and frequency spaces
respectively. Basically we have no a priori control on these
parameters, which is the main source of the difficulty of
establishing useful properties for almost periodic modulo symmetry
solutions. However, under the spherical symmetry assumption, one is
allowed to fix the center of mass, thus leaving only one parameter
$N(t)$ which can still vary arbitrarily. This case turns out to be
treatable in high dimensions $d\ge 4$. Here is the main theorem of
this paper:

\begin{thm}\label{main}
Let $d\ge 4$. Let $u$ be a maximal-lifespan solution on $I$ and is
spherically symmetric. Suppose $u$ is almost periodic modulo scaling
on $I$. Then there exists $\eps=\eps(d)<\frac 4d$ such that
\begin{align}\label{add_reg}
u(t)\in H_x^{1+\eps}, \ \forall \ t\in I.
\end{align}
Moreover, the kinetic energy of the solution is localized uniformly
in time: for any $\eta>0$, there exists $C(\eta)$ such that for
any $t\in I$
\begin{align}\label{kinetic_loc}
\int_{|x|\ge C(\eta)}|\nabla u(t,x)|^2 dx \le \eta.
\end{align}
Here, $\eps$ only depends on the dimension $d$, while $C(\eta)$
depends also on the solution $u$.
\end{thm}

A few remarks are in order.

\begin{rem}
This result seems a bit surprising in view of the fact that the
scaling parameter $N(t)$ can vary arbitrarily and the solution is
only assumed to be in the scale invariant space $L_x^2$. On the
other hand, Theorem \ref{main} bears similarities with previous
works \cite{klvz,lz:2d} and \cite{kvz:blowup, ktv:2d} where they
were able to deal with dimensions two and higher. However in
\cite{klvz, lz:2d}, the solution is assumed to have $H_x^1$
regularity and this latter fact allows one to treat solutions being
almost periodic modulo scaling in only one time direction. In
\cite{kvz:blowup, ktv:2d}, the additional regularity is only
established for three typical solutions known as \emph{three
enemies}. Namely, these
are \emph{almost periodic modulo scaling} solutions with a priori control on $N(t)$: \\
\indent a. The self-similar solution. This solution is defined on maximal time interval $(0,\infty)$ and $N(t)=t^{-\frac 12}$ for any $t\in (0,\infty)$;

b. The soliton-like solution. This solution is global and $N(t)=1$;

c. The high to low cascade. This solution is also global with $N(t)$ satisfying: $N(t)\le 1$, $\liminf_{t\to \pm\infty} N(t)=0$.

\noindent On the other hand, the technique in this paper allows us to deal with all enemies \emph{with no a priori assumption on $N(t)$} in dimensions $d\ge4$.
\end{rem}

%\begin{rem}
%This result bares similarities with the previous work
%\cite{kvz:blowup, ktv:2d}. However, there the additional regularity
%is only established for three typical solutions referred to as \emph{three enemies}. Namely, the
%\emph{almost periodic modulo scaling} solutions with a priori
%control on $N(t)$: the self similar solution with $N(t)=t^{-\frac
%12}$; the soliton-like solution with $N(t)=1$ and the high to low
%cascade with $N(t)\le 1$ and $\liminf_{t\to \pm\infty} N(t)=0$. On the
%other hand, the technique in this paper allows us to deal with all
%the enemies \emph{with no a priori assumption on $N(t)$}.
%\end{rem}

\begin{rem}
The dependence on the dimension comes from the fact that in dimension $d\ge 4$, the nonlinearity $|u|^{\frac 4d}u$ can be put in Lebesgue
space $L_x^p(\R^d)$ for some $p\ge 1$ only knowing that $u\in L_x^2(\R^d)$. This property is not available in low dimensions $d=2,3$. So in these dimensions, it is still open proving the additional regularity for solutions other than the three enemies. \end{rem}

\begin{rem}
Besides the spherical symmetry, we can also consider other
symmetries that can freeze the center of mass at the origin. For
example, one can consider the splitting spherical symmetry
introduced in \cite{lz:2d}. In \cite{lz:six}, we select  the six
dimensions as a sample case to show how the technique can be
extended to deal with the solution with splitting spherical symmetry
and is almost periodic modulo scaling. There the main difficulty
comes from the fact that the waves can propagate anisotropically
along splitting subspaces. As shown in the proof of Proposition
\ref{fre_decay} and Proposition \ref{spa_decay}, the spherical
symmetry is mainly used to treat the part where the plane waves
travel away from the origin. For this part, one uses the weighted
Strichartz estimate for radial functions to get the decay. In the
splittingly spherical symmetric case, we develop tools such as
weighted Strichartz estimate (see \cite{lz:2d}) for splittingly
spherical symmetric functions to make use of the decay property.
\end{rem}

\begin{rem} \label{rem412}
To prove Theorem \ref{main} we need to control the parts of the
solution both near the spatial origin and away from it. To control
the part away from the origin, we use the techniques from
\cite{klvz} where we need the radial assumption on the solution. To
control the part near the origin,  we introduce a novel \emph{local
iteration} scheme which actually does not need the radial assumption
provided we already have the control on the piece away from the
origin.
%The only part where the radial assumption comes is in the first part of the proof, where we need to control the piece away
%from the origin.
 We should also stress that our proof uses the almost periodicity in a very light way. Instead of assuming the
solution is almost periodic modulo scaling on the whole time interval, one could assume the following
\emph{sequential almost periodicity}: there exist $t_n^{+} \to \sup I$, $t_n^{-} \to \inf I$ and scaling parameters
$N(t_n^{+})$, $N(t_n^{-})$, such that both of the sets
$$\{ N(t_n^{+})^{-\frac d2} u(t_n^{+}, \cdot/N(t_n^{+})) \},\, \{ N(t_n^{-})^{-\frac d2} u(t_n^{-}, \cdot/N(t_n^{-}))\}$$
are precompact in $L_x^2(\R^d)$.
\end{rem}
\subsection{Applications of Theorem \ref{main}}

The applications of Theorem \ref{main} are related to the scattering conjecture and the rigidity conjecture which we now explain. In the defocusing case, the scattering conjecture says that all solutions with finite mass exist globally and scatter in both time directions. In the focusing case, besides scattering solutions, there exist finite time blowup solutions as shown in \cite{glassey} and the solitary wave solutions of the form $e^{it}R(x)$.  Here $R$ solves the elliptic equation
\begin{align*}
\Delta R-R+|R|^{\frac 4d}R=0.
\end{align*}
There are infinitely many solutions to this equation, but only one positive solution which
is spherically symmetric (up to translations) and whose mass is minimal among all these $R's$. This solution is usually called the

\begin{defn}[Ground state \cite{blions, kwong}] \label{ground} The ground state $Q$ refers to the unique positive radial Schwartz solution to the elliptic equation
\begin{align*}
\Delta Q-Q+|Q|^{\frac 4d}Q=0.
\end{align*}
\end{defn}

It is believed that the mass of $Q$ serves as the minimal mass among all the nonscattering solutions in the focusing case. To summarize, we have

\begin{conj}[Scattering conjecture] Let $u_0\in L_x^2(\R^d)$. In the focusing case, we also assume $M(u_0)<M(Q)$. Then the corresponding solution to \eqref{nls} exists globally and scatters in both time directions.
\end{conj}

This conjecture has been proved in dimensions $d\ge 2$ when the
initial data is spherically symmetric, see \cite{ktv:2d},
\cite{kvz:blowup}.\footnote{In the defocusing case and  $d\ge 3$,
one can take advantage of Morawetz estimate to prove the additional
regularity, see \cite{tvz:hd} for more details.} We now give a high
level overview of the proof which is based on a contradiction
argument. Assume the scattering conjecture is not true, one can then
use concentration compactness tools to obtain minimal mass
non-scattering\footnote{Here by "non-scattering", we mean that the
$L_{t,x}^{\frac{2(d+2)}d}$ norm of the solution is infinite.
Obviously, a "non-scattering" solution may blow up at finite time or
exist globally with infinite $L_{t,x}^{\frac{2(d+2)}d}$ norm. }
solutions which are almost periodic modulo scaling (due to the
spherical symmetry) with scaling parameter $N(t)$. To obtain better
control of $N(t)$, another limiting procedure is performed to reduce
the consideration to three typical solutions alluded as to "three
enemies". To kill three enemies and thereby obtaining the
contradiction, one can use the information of $N(t)$ to obtain
additional regularity of these solutions which together with a
truncated virial argument establishes the claim.

Thanks to Theorem \ref{main}, we can simplify the argument by
directly working with all enemies whose scaling parameter $N(t)$ can
vary arbitrarily in dimensions $d\ge 4$. In other words, the
limiting procedure of picking three enemies is not needed here. We
record the result as:

\begin{cor}[Scattering in dimension $d\ge 4$ with spherical symmetry]\label{scatter}
Let $d\ge 4$. Let $u_0\in L_x^2(\R^d)$ be spherically symmetric. In the focusing case, we assume $M(u_0)<M(Q)$. Then the solution to \eqref{nls} with this initial data exists globally and satisfies
\begin{align*}
\|u\|_{L_{t,x}^{\frac{2(d+2)}d}(\R\times\R^d)}\le  C(\|u_0\|_{L_x^2}).
\end{align*}
\end{cor}

\vspace{0.2cm}

We turn now to the rigidity conjecture.

In the focusing case, a main issue is to understand the large time
behavior of non-scattering solutions. This problem has only been
addressed in the case when the mass of $u$ is equal to or slightly
bigger than that of the ground state, see \cite{merle-raph},
\cite{merle_duke}, \cite{klvz}, \cite{lz:2d} and the references
therein. In this paper, we are primarily concerned with the case
when the solution has the ground state mass. Our main focus is to
characterize and classify all such solutions. At the level of ground
state mass, there are two explicit examples of non-scattering
solutions: the solitary wave $SW$ which exists globally and the
pseudo-conformal ground state $Pc(Q)$ which blows up at $t=0$:

\begin{align*}
SW&=e^{it}Q(x),\\
Pc(Q)&=|t|^{-\frac d2}e^{\frac{i|x|^2-4}{4t}}Q(\frac x{t}).
\end{align*}
It is conjectured that, up to symmetries, these are the only two
threshold solutions for scattering at the level of minimal mass.
Associated with this is the following rigidity conjecture which
identifies all solutions with ground state mass as either $SW$ or
$Pc(Q)$ if they do not scatter. Since both mass and the equation are
invariant under a couple of symmetries, the coincidence of the
solutions with the examples only hold modulo these symmetries.
Specifically, the symmetries are: translation,
 phase rotation, scaling and the Galilean boost.

\begin{conj}[Rigidity conjecture at the ground state mass]
Let $u_0\in L_x^2(\R^d)$ satisfy $M(u_0)=M(Q)$. Then only the
following cases can occur

1. The solution $u$ blows up at finite time, then in this case $u$
must coincide with $Pc(Q)$ up to symmetries of the equation.

2. The solution $u$ is a global solution. Then in this case, $u$
either scatters in both time directions or $u$ must coincide with
$SW$ up to symmetries of the equation.
\end{conj}

In \cite{merle_duke}, Merle considered the first part of the
conjecture, where he identified all finite time blowup solutions as
$Pc(Q)$ under an additional $H_x^1$ assumption on the initial data. See also
\cite{weinstein:charact} for the preliminary result due to Weinstein
and \cite{hmidi-keraani} for a simplified proof of Merle's argument
due to Hmidi-Keraani. By Merle's result and pseudoconformal
transformation, the second part of the conjecture, which
characterizes all global solutions with ground state mass, still
holds if we make the strong assumption that the initial data $u_0\in
\Sigma=\{f\in H_x^1,xf\in L_x^2\}$. Finally it is worthwhile
noticing that Merle's argument works for all dimensions without any
symmetry assumption on the initial data.

Without the $\Sigma$ assumption on the initial data, it is not clear
at all how to deal with the case when $u_0$ is merely in $L_x^2$ and
the corresponding solution is global. Recently in \cite{klvz} and
\cite{lz:2d}, we proved the second part of the conjecture when the
initial data $u_0\in H_x^1(\R^d)$, $d\ge 2$ and is spherically
symmetric. In dimension $d\ge 4$, the results hold even under a
weaker symmetry assumption, namely, the initial data is only
required to be splitting-spherical symmetric (see \cite{lz:2d} for
more details).

As stated, all the results concerning the rigidity conjecture
require the $H_x^1$ regularity on the initial data since it is the
minimal regularity to define the energy and to carry out the
spectral analysis. Here the energy refers to
\begin{align*}
\mbox{Energy:}\ \ E(u(t))=\frac 12\|\nabla u(t)\|_{L_x^2}^2-\frac
d{2(d+2)} \|u(t)\|_{L_{x}^{\frac{2(d+2)}d}}^{\frac{2(d+2)}d}=E(u_0).
\end{align*}
To prove the rigidity results for pure $L_x^2$ solutions, a reasonable strategy is to upgrade the regularity of the solution to $H_x^1$ or better by taking advantage of certain compactness properties of the solutions. This is where Theorem \ref{main} has to be used. We can then use known $H_x^1$ results to classify these solutions. Therefore as a direct consequence of Theorem \ref{main}, we have

\begin{thm}[Rigidity for two-way non-scattering solutions with ground state mass]\label{rigid}
Let $d\ge 4$. Let $u_0\in L_x^2(\R^d)$ be spherically symmetric and $M(u_0)=M(Q)$. Let $u$ be the maximal lifespan solution on $I$ which does not scatter on both sides:
\begin{align*}
\|u\|_{L_{t,x}^{\frac{2(d+2)}d}([t_0,\sup I)\times\R^d)}=\|u\|_{L_{t,x}^{\frac{2(d+2)}d}
((\inf I, t_0]\times\R^d)}=\infty, \ t_0\in I.
\end{align*}
Then $I=\R$ and $u=e^{it}Q$ up to phase rotation and scaling.
\end{thm}

For technical reasons, we need to impose the condition that the solution does not scatter in both time directions. It is an interesting problem to extend our techniques to the case  when the solution scatters only in one time direction, but does not scatter in the other.

We will give the proof of these two results in section \ref{apply}. Now we briefly sketch the proof of Theorem \ref{main}.

\subsection{Main idea of the proof of Theorem \ref{main}: a local iteration scheme}

We will work with each single dyadic
frequency of $u$:
\begin{align*}
 \|P_N u(t)\|_{L_x^2}.
\end{align*}
The decay in $N$ will correspond to the regularity of the solution. First we observe that
when restricted to the region away from the origin, the argument in \cite{klvz} gives us
\begin{align} \label{eq_tmp_542}
\|\phi_{>1}P_N u(t)\|_{L_x^2}\lsm N^{-1-\eps}
\end{align}
with a uniform in time bound. Here $\phi_{>1}$ is a smooth cut-off function supported in the
region $|x|>1$. This reduces matters to estimating the part of
the solution near the spatial origin, i.e. $\|\phi_{\le 1}P_N
u(t)\|_{L_x^2}$. This piece is trivially bounded by
\begin{align*}
A_N=\|P_N u\|_{S([t,t+\frac 1{\sqrt N}])},
\end{align*}
i.e.
the Strichartz norm
of $P_Nu$ on a local time interval $[t, t+\frac 1{\sqrt N}]$.
It turns out, after some technical manipulations, that this latter quantity is better suited for iteration and
bootstrapping. Indeed we shall establish recurrent relations for $A_N$ and we will iterate
our estimates only finitely many (but sufficiently many)
steps. The crucial point is that during the iteration process, we
shall never need more than the information of the solution on a unit time
interval $[t,t+1]$. Therefore we do not need to use the full control
of $N(t)$. We remark that although as a sacrifice the $H_x^{1+\eps}$
norm of $u(t)$ depends on $t$, this information combined with the kinetic energy localization in Section \ref{apply} suffice to prove Corollary \ref{scatter} and Theorem \ref{rigid}.

\subsection*{Acknowledgements}
Both authors were supported in part by the National Science Foundation under
agreement No. DMS-0635607 and a start-up funding from Mathematics
Department of University of Iowa. X.~Zhang was also supported by NSF
grant No.~10601060 and project 973 in China. D. Li was also supported by
NSF grant No.~0908032 and the old gold summer fellowship from University of Iowa.

%%%%%%%%%%%%%%%%%%%%%%%%%%%%%%%%%%%%%%%%%%%%%%%%%%%%%%%%%%%%%%%%%%%%%%%%%%%%%%%%%%%%%%%%%%%
%
%
%                                   Section
%
%
%%%%%%%%%%%%%%%%%%%%%%%%%%%%%%%%%%%%%%%%%%%%%%%%%%%%%%%%%%%%%%%%%%%%%%%%%%%%%%%%%%%%%%%%%%%

\section{Preliminaries}

\subsection{Some notations}
We write $X \lesssim Y$ or $Y \gtrsim X$ to indicate $X \leq CY$ for some constant $C>0$.  We use $O(Y)$ to denote any quantity $X$
such that $|X| \lesssim Y$.  We use the notation $X \sim Y$ whenever $X \lesssim Y \lesssim X$.  The fact that these constants
depend upon the dimension $d$ will be suppressed.  If $C$ depends upon some additional parameters, we will indicate this with
subscripts; for example, $X \lesssim_u Y$ denotes the assertion that $X \leq C_u Y$ for some $C_u$ depending on $u$. Sometimes
when the context is clear, we will suppress the dependence on $u$ and write $X \lesssim_u Y$ as $X \lesssim Y$.
We will write $C=C(Y_1, \cdots, Y_n)$ to stress that the constant $C$ depends on quantities $Y_1$, $\cdots$, $Y_n$.
We denote by $X\pm$ any quantity of the form $X\pm \epsilon$ for any $\epsilon>0$.

We use the `Japanese bracket' convention $\langle x \rangle := (1 +|x|^2)^{1/2}$.

We write $L^q_t L^r_{x}$ to denote the Banach space with norm
$$ \| u \|_{L^q_t L^r_x(\R \times \R^d)} := \Bigl(\int_\R \Bigl(\int_{\R^d} |u(t,x)|^r\ dx\Bigr)^{q/r}\ dt\Bigr)^{1/q},$$
with the usual modifications when $q$ or $r$ are equal to infinity, or when the domain $\R \times \R^d$ is replaced by a smaller
region of spacetime such as $I \times \R^d$.  When $q=r$ we abbreviate $L^q_t L^q_x$ as $L^q_{t,x}$.

Throughout this paper, we will use $\phi\in C^\infty(\R^d)$ be a
radial bump function supported in the ball $\{ x \in \R^d: |x| \leq
\frac{25} {24} \}$ and equal to one on the ball $\{ x \in \R^d: |x|
\leq 1 \}$.  For any constant $C>0$, we denote $\phi_{\le C}(x):=
\phi \bigl( \tfrac{x}{C}\bigr)$ and $\phi_{> C}:=1-\phi_{\le C}$.

\subsection{Basic harmonic analysis}\label{ss:basic}

For each number $N > 0$, we define the Fourier multipliers
\begin{align*}
\widehat{P_{\leq N} f}(\xi) &:= \phi_{\leq N}(\xi) \hat f(\xi)\\
\widehat{P_{> N} f}(\xi) &:= \phi_{> N}(\xi) \hat f(\xi)\\
\widehat{P_N f}(\xi) &:= (\phi_{\leq N} - \phi_{\leq N/2})(\xi) \hat
f(\xi)
\end{align*}
and similarly $P_{<N}$ and $P_{\geq N}$.  We also define
$$ P_{M < \cdot \leq N} := P_{\leq N} - P_{\leq M} = \sum_{M < N' \leq N} P_{N'}$$
whenever $M < N$.  We will usually use these multipliers when $M$ and $N$ are \emph{dyadic numbers} (that is, of the form $2^n$
for some integer $n$); in particular, all summations over $N$ or $M$ are understood to be over dyadic numbers.  Nevertheless, it
will occasionally be convenient to allow $M$ and $N$ to not be a power of $2$.  As $P_N$ is not truly a projection, $P_N^2\neq P_N$,
we will occasionally need to use fattened Littlewood-Paley operators:
\begin{equation}\label{PMtilde}
\tilde P_N := P_{N/2} + P_N +P_{2N}.
\end{equation}
These obey $P_N \tilde P_N = \tilde P_N P_N= P_N$.

Like all Fourier multipliers, the Littlewood-Paley operators commute with the propagator $e^{it\Delta}$, as well as with
differential operators such as $i\partial_t + \Delta$. We will use basic properties of these operators many times,
including

\begin{lem}[Bernstein estimates]\label{Bernstein}
 For $1 \leq p \leq q \leq \infty$,
\begin{align*}
\bigl\| |\nabla|^{\pm s} P_N f\bigr\|_{L^p_x(\R^d)} &\sim N^{\pm s} \| P_N f \|_{L^p_x(\R^d)},\\
\|P_{\leq N} f\|_{L^q_x(\R^d)} &\lesssim N^{\frac{d}{p}-\frac{d}{q}} \|P_{\leq N} f\|_{L^p_x(\R^d)},\\
\|P_N f\|_{L^q_x(\R^d)} &\lesssim N^{\frac{d}{p}-\frac{d}{q}} \| P_N f\|_{L^p_x(\R^d)}.
\end{align*}
\end{lem}

While it is true that spatial cutoffs do not commute with Littlewood-Paley operators, we still have the following:

\begin{lem}[Mismatch estimates in real space]\label{L:mismatch_real}
Let $R,N>0$.  Then
\begin{align*}
\bigl\| \phi_{> R} \nabla P_{\le N} \phi_{\le\frac R2} f \bigr\|_{L_x^p(\R^d)}  &
\lsm_m N^{1-m} R^{-m} \|f\|_{L_x^p(\R^d)} \\
\bigl\| \phi_{> R}  P_{\leq N} \phi_{\le\frac R2} f \bigr\|_{L_x^p(\R^d)}
&\lsm_m N^{-m} R^{-m} \|f\|_{L_x^p(\R^d)}
\end{align*}
for any $1\le p\le \infty$ and $m\geq 0$.
\end{lem}

\begin{proof}
We will only prove the first inequality; the second follows similarly.

It is not hard to obtain kernel estimates for the operator $\phi_{>
R}\nabla P_{\le N}\phi_{\le\frac R2}$. Indeed, an exercise in
non-stationary phase shows
\begin{align*}
\bigl|\phi_{> R}\nabla P_{\le N}\phi_{\le\frac R2}(x,y)\bigr|
\lesssim N^{d+1-2k} |x-y|^{-2k}\phi_{|x-y|>\frac R2}
\end{align*}
for any $k\geq 0$.  An application of Young's inequality yields the claim.
\end{proof}

Similar estimates hold when the roles of the frequency and physical spaces are interchanged.  The proof is easiest when
working on $L_x^2$, which is the case we will need; nevertheless, the following statement holds on $L_x^p$ for any $1\leq p\leq \infty$.

\begin{lem}[Mismatch estimates in frequency space]\label{L:mismatch_fre}
For $R>0$ and $N,M>0$ such that $\max\{N,M\}\geq 4\min\{N,M\}$,
\begin{align*}
\bigl\|  P_N \phi_{\le{R}} P_M f \bigr\|_{\ltrd} &\lsm_m \max\{N,M\}^{-m} R^{-m} \|f\|_{\ltrd} \\
\bigl\|  P_N \phi_{\le {R}} \nabla P_M f \bigr\|_{\ltrd} &\lsm_m M
\max\{N,M\}^{-m} R^{-m} \|f\|_{\ltrd}.
\end{align*}
for any $m\geq 0$.  The same estimates hold if we replace $\phi_{\le
R}$ by $\phi_{>R}$.
\end{lem}

\begin{proof}
The first claim follows from Plancherel's Theorem, Lemma~\ref{L:mismatch_real} and its adjoint.  To obtain the second claim from this, we write
$$
P_N \phi_{\le {R}} \nabla P_M = P_N \phi_{\le {R}} P_M \nabla \tilde
P_M
$$
and note that $\|\nabla \tilde P_M\|_{L_x^2\to L_x^2}\lesssim M$.
\end{proof}

\subsection{Some analysis tools}

We will need the following radial Sobolev embedding to exploit the
decay property of a radial function. For the proof and the more
complete version, see \cite{tvz:hd}.

\begin{lem}[Radial Sobolev embedding, \cite{tvz:hd}]\label{L:radial_embed}
Let dimension $d\ge 2$. Let $s>0$, $\alpha>0$, $1<p,q<\infty$ obeys
the scaling restriction: $\alpha+s=d(\frac 1q-\frac 1p)$. Then the
following holds:
\begin{align*}
\||x|^{\alpha} f\|_{L^p(\R^d)}\lsm \||\nabla|^s f\|_{L^q(\R^d)},
\end{align*}
where the implicit constant depends on $s,\alpha,p,q$. When
$p=\infty$, we have
\begin{align*}
\||x|^{\frac{d-1}2}P_N f\|_{L^\infty(\R^d)}\lsm N^{\frac 12}\|P_N
f\|_{\ltrd}.
\end{align*}
\end{lem}

We will need the following fractional chain rule lemma.
\begin{lem}[Fractional chain rule for a $C^1$ function, \cite{chris:weinstein}\cite{staf}\cite{taylor}]\label{lem_chain}
Let $G\in C^1(\mathbb C)$, $\sigma \in (0,1)$, and $1<r,r_1,r_2<\infty$ such that $\frac 1r=\frac 1{r_1}+\frac 1{r_2}$.
Then we have
\begin{align*}
\||\nabla|^{\sigma}G(u)\|_r\lesssim \|G'(u)\|_{r_1}\||\nabla|^{\sigma}u\|_{r_2}.
\end{align*}
\end{lem}
\begin{proof}
See \cite{chris:weinstein}, \cite{staf} and \cite{taylor}.
\end{proof}

We also need the following lemma from \cite{kvz:blowup}.
\begin{lem}\label{special}
Let $0<s<1+\frac 4d$ and $F(u)=|u|^{\frac 4d}u$. Then
\begin{align*}
\||\nabla|^s F(u)\|_{L_x^{\frac{2(d+2)}{d+4}}}\lsm\||\nabla|^s
u\|_{L_x^{\frac{2(d+2)}d}} \|u\|_{L_x^{\frac{2(d+2)}d}}^{\frac 4d}.
\end{align*}
\end{lem}

We will need the following sharp Gagliardo-Nirenberg inequality
\begin{lem}[\cite{W1}]\label{L:sharp_gn} Let $Q$ be the ground state in the
Definition \ref{ground}. Then for any $f\in H_x^1(\R^d)$, we have
\begin{align}\label{sharp-gn}
\|f\|_{L_x^{\frac{2(d+2)}d}}^{\frac{2(d+2)}d}\le
\frac{d+2}d\biggl(\frac{M(f)}{M(Q)}\biggr)^{\frac 2d}\|\nabla
f\|_{L_x^2}^2.
\end{align}
The equality holds only and if only
\begin{align}
f=ce^{i\theta}\lambda^{\frac d2} Q(\lambda(x-x_0))
\end{align}
for $(c,\theta,\lambda)\in (\R^+,\R,\R^+)$.
\end{lem}

\subsection{Strichartz estimates}

The free Schr\"odinger flow has the explicit expression:
\begin{align*}
e^{it\Delta } f(x)=\frac 1{(4\pi t)^{d/2}}\int_{\R^d}
e^{i|x-y|^2/4t}f(y)dy,
\end{align*}
%from which we also have the dispersive estimate
%\begin{align*}
%\|e^{it\Delta}f\|_{L^{\infty}(\R^d)}\lsm \frac
%1{|t|^{d/2}}\|f\|_{L^1(\R^d)}.
%\end{align*}
from which we can derive the kernel estimate of the frequency
localized propagator. We record the following

\begin{lem}[Kernel estimates, \cite{ktv:2d,kvz:blowup}]\label{L:kernel}
For any $m\ge 0$, we have
\begin{align*}
|(\pn e^{it\Delta}(x,y)|\lsm_m
\begin{cases}
|t|^{-d/2},&: |x-y|\sim Nt;\\
\frac{N^d}{|N^2 t|^m\langle N|x-y|\rangle^m}&: \mbox{otherwise}
\end{cases}
\end{align*}
for $|t|\ge N^{-2}$ and
\begin{align*}
|(\pn e^{it\Delta})(x,y)|\lsm_m N^d\langle N|x-y|\rangle^{-m}
\end{align*}
for $|t|\le N^{-2}$.
\end{lem}

We will frequently use the standard Strichartz estimate. Let $d\ge 3$. Let $I$ be a time interval. We define the Strichartz space on $I$:
$$
S(I)=L_t^{\infty}L_x^2(I\times\R^d)\cap L_t^2L_x^{\frac{2d}{d-2}}(I\times\R^d).
$$
We also define $N(I)$ to be $L_t^1L_x^2(I\times\R^d)+ L_t^2 L_x^{\frac{2d}{d+2}}(I\times\R^d)$. Then the standard Strichartz estimate reads

\begin{lem}[Strichartz]\label{L:strichartz}
Let $d\ge 3$. Let $I$ be an interval, $t_0 \in I$, and let $u_0 \in L^2_x(\R^d)$
and $f \in N(I)$.  Then, the
function $u$ defined by
$$
u(t) := e^{i(t-t_0)\Delta} u_0 - i \int_{t_0}^t e^{i(t-t')\Delta} f(t')\ dt'
$$
obeys the estimate
$$
\|u \|_{S(I)} \lesssim \| u_0 \|_{L^2_x} + \|f\|_{N(I)},
$$
where all spacetime norms are over $I\times\R^d$.
\end{lem}

\begin{proof}
See, for example, \cite{gv:strichartz, strichartz}.  For the
endpoint see \cite{tao:keel}.
\end{proof}

We will also need a weighted Strichartz estimate, which exploits heavily the spherical symmetry in order to obtain spatial decay.

\begin{lem}[Weighted Strichartz, \cite{ktv:2d, kvz:blowup}]\label{L:wes} Let $I$ be an interval, $t_0 \in I$, and let
$F:I\times\R^d\to \C$ be spherically symmetric.  Then,
$$ \biggl\|\int_{t_0}^t e^{i(t-t')\Delta} F(t')\, dt' \biggr\|_{L_x^2}
\lesssim \bigl\||x|^{-\frac{2(d-1)}q}F \bigr\|_{L_t^{\frac{q}{q-1}}L_x^{\frac{2q}{q+4}}(I \times \R^d)}
$$
for all $4\leq q\leq \infty$.
\end{lem}

\subsection{The in-out decomposition}
We will need an incoming/outgoing decomposition; we will use the one developed in \cite{ktv:2d, kvz:blowup}.
As there, we define operators $P^{\pm}$ by
\begin{align*}
[P^{\pm} f](r) :=\tfrac12 f(r)\pm \tfrac{i}{\pi} \int_0^\infty \frac{r^{2-d}\,f(\rho)\,\rho^{d-1}\,d\rho}{r^2-\rho^2},
\end{align*}
where the radial function $f: \R^d\to \C$ is written as a function of radius only.
We will refer to $P^+$ is the projection onto outgoing spherical waves; however, it is not a true projection as it is neither idempotent
nor self-adjoint.  Similarly, $P^-$ plays the role of a projection onto incoming spherical waves; its kernel is the complex
conjugate of the kernel of $P^+$ as required by time-reversal symmetry.

For $N>0$ let $P_N^{\pm}$ denote the product $P^{\pm}P_N$ where $P_N$ is the Littlewood-Paley projection.
We record the following properties of $P^{\pm}$ from \cite{ktv:2d, kvz:blowup}:

\begin{prop}[Properties of $P^\pm$, \cite{ktv:2d, kvz:blowup}]\label{P:P properties}\leavevmode
\begin{SL}
\item $P^+ + P^- $ represents the projection from $L^2$ onto
$L^2_{\text{rad}}$.

\item Fix $N>0$.  Then
$$
\bigl\| \chi_{\gtrsim\frac 1N} P^{\pm}_{\geq N} f\bigr\|_{L^2(\R^d)}
\lesssim \bigl\| f \bigr\|_{L^2(\R^d)}
$$
with an $N$-independent constant.

\item If the dimension $d=2$, then the $P^{\pm}$ are bounded on $L^2(\R^2)$.
\item For $|x|\gtrsim N^{-1}$ and $t\gtrsim N^{-2}$, the
integral kernel obeys
\begin{equation*}
\bigl| [P^\pm_N e^{\mp it\Delta}](x,y) \bigr| \lesssim \begin{cases}
    (|x||y|)^{-\frac {d-1}2}|t|^{-\frac 12}  &: \  |y|-|x|\sim  Nt \\[1ex]
     \frac{N^d}{(N|x|)^{\frac{d-1}2}\langle N|y|\rangle^{\frac{d-1}2}}
     \bigl\langle N^2t + N|x| - N|y| \bigr\rangle^{-m}
            &: \  \text{otherwise}\end{cases}
\end{equation*}
for all $m\geq 0$.

\item For $|x|\gtrsim N^{-1}$ and $|t|\lesssim N^{-2}$, the
integral kernel obeys
\begin{equation*}
\bigl| [P^\pm_N e^{\mp it\Delta}](x,y) \bigr|
    \lesssim   \frac{N^d}{(N|x|)^{\frac{d-1}2}\langle N|y|\rangle^{\frac{d-1}2}}
     \bigl\langle N|x| - N|y| \bigr\rangle^{-m}
\end{equation*}
for any $m\geq 0$.
\end{SL}
\end{prop}

\section{Theorem \ref{main} implies Corollary
\ref{scatter} and Theorem \ref{rigid}}\label{apply}

In this section, we assume Theorem \ref{main} holds momentarily and
prove the scattering and the rigidity result Corollary
\ref{scatter}, Theorem \ref{rigid}.

\textbf{Proof of Corollary \ref{scatter}:}
\begin{proof} Suppose by contradiction that Corollary
\ref{scatter} does not hold. Then
there exist minimal mass $M_c$ for which $M_c<\infty$ in the defocusing case, $M_c<M(Q)$
in the focusing case and maximal-lifespan solution $u(t,x)$ on
$I=(-T_*, T^*)$ such that

1. $u$ is spherically symmetric and  $M(u)=M_c$;

2. $u$ is almost periodic modulo scaling on $I$.

See for instance \cite{compact} for this part of the argument which is by now standard.
Applying Theorem \ref{main}, we know that $u\in H_x^{1+\eps}$. We now detail the rest of the argument in the focusing case,
since the defocusing case is even simpler. By the sharp Gagliardo-Nirenberg inequality and the fact that
$M(u)<M(Q)$ we have
\begin{align*}
\|u(t)\|_{H_x^1}\lsm_{M(u)} 1.
\end{align*}
From this and the standard local theory in $H_x^1$ we know that $u$
exists globally, ie: $T_*=T^*=\infty$. In this situation, the
contradiction will come from the truncated virial and the kinetic
energy localization as we explain now. Let $\phi_{\le R}$ be the
smooth cutoff function, we define the truncated virial as
\begin{align*}
V_R(t)=\int \phi_{\le R}(x) |x|^2 |u(t,x)|^2 dx.
\end{align*}
Obviously
\begin{align}\label{vr_bdd}
 V_R(t)\lsm R^2, \ \forall\ t\in \R.
\end{align}
On the other hand, we compute the second derivative of virial with
respect to $t$, this gives
\begin{align}\label{vir}
\partial_{tt}& V_R(t)=8 E(u)+\\
&O\biggl(\int_{|x|>R}|\nabla u(t,x)|^2 +|u(t,x)|^{\frac {2(d+2)}d}+\frac
1{R^2}\int_{|x|>R}|u(t,x)|^2 dx\biggr).\label{error}
\end{align}
Since $M(u)<M(Q)$ and $u\in H_x^1$, from sharp Gagliardo-Nirenberg
inequality \eqref{sharp-gn} we have
\begin{align*}
E(u)>0.
\end{align*}
Now we can use the kinetic energy localization \eqref{kinetic_loc}
and Gagliardo-Nirenberg inequality to control the error term
\eqref{error} and finally get
\begin{align*}
\partial_{tt} V_R(t)\ge 4 E(u)>0
\end{align*}
by taking $R$ sufficiently large.  This obviously contradicts
\eqref{vr_bdd}. The proof of Corollary \ref{scatter} is finished.
\end{proof}

\textbf{The proof of Theorem \ref{rigid}}

\begin{proof} Let $d\ge 4$, let $u$ be the solution of \eqref{nls}
satisfying the following

1. $M(u)=M(Q)$, $u$ is spherically symmetric.

2. $u$ does not scatter in both time directions.

By \cite{kvz:blowup} or Corollary \ref{scatter}, $M(Q)$ is the minimal mass, the
compactness argument in \cite{compact,keraani,BegoutVargas} shows
that $u$ is \emph{almost periodic modulo scaling} in both time
directions. Now we can apply Theorem \ref{main} to deduce that $u\in
H_x^1$. Since from Merle's result, the only finite time blowup
solution must be $Pc(Q)$ up to symmetries and $Pc(Q)$ scatters in
one time direction, we know from condition 2 that $u$ must be a
global solution.

From \eqref{sharp-gn}, this global solution $u$ satisfies $E(u)\ge
0$. Moreover, the same virial argument as in the proof of Corollary
\ref{scatter} precludes the case $E(u)>0$, we therefore obtain $E(u)=0$.
From here the coincidence of the solution with solitary wave follows
immediately, again by the sharp Gagliardo-Nirenberg inequality.

\end{proof}

\section{The proof of Theorem \ref{main}}

  The proof of Theorem \ref{main} proceeds in two steps. In
the first step, we prove that away from the origin, the solution has $H_x^{1+\eps}$ regularity. Moreover, a
similar (but more refined) argument establishes the spatial decay estimate. These two
pieces together suffice for us to establish the kinetic localization
estimate. However, in this step, the total kinetic energy does not
need to be finite.

In the second step, we prove the total kinetic energy is actually finite by controlling the piece near the spatial origin.
Thanks to the first step, we only need to consider a single
frequency $P_N u$ with spatial cutoff $\phi_{\le 1}$. We can bound this quantity by the Strichartz norm of $P_N u$ on a short time
interval $[t,t+\frac 1{\sqrt N}]$. We then establish a recurrent
relation for this local Strichartz norm.
Iterating the estimates finitely many times then yields the desired bound. More details are given below.

Before proceeding, we remark that in all of the arguments that follow,
the only property we use for an \emph{almost periodic modulo scaling} solution is that it satisfies
the improved Duhamel formula. This was first derived
in \cite{compact}.

\begin{prop}[Improved Duhamel formula, \cite{compact}]\label{impr_duham}
Let $u$ be the solution of \eqref{nls} and is almost periodic modulo
scaling on the time interval $I$. Then we have the following
\begin{align}\label{improve}
u(t)&=\lim_{T\to\inf I} -i\int_{T}^t e^{i(t-\tau)\Delta}
F(u(\tau))d\tau\\
&=\lim_{T\to\sup I} i\int_t^{T} e^{i(t-\tau)\Delta}
F(u(\tau))d\tau.\notag
\end{align}

Here the limit is in weak $L_x^2$ sense.
\end{prop}

\begin{rem}
As was already mentioned in Remark \ref{rem412}, we actually only need the
\emph{sequential almost periodicity} of the solution for the later proof to work.
This would imply the following sequence version of improved Duhamel formula:
\begin{align*}%\label{improve1}
u(t)&=\lim_{n \to \infty} -i\int_{T_n^{-} }^t e^{i(t-\tau)\Delta}
F(u(\tau))d\tau\\
&=\lim_{n \to \infty} i\int_t^{T_n^{+}} e^{i(t-\tau)\Delta}
F(u(\tau))d\tau.\notag
\end{align*}
Here again the limit is in weak $L_x^2$ sense.

\end{rem}
\noindent
In what follows, we shall only assume $u$ satisfy
\begin{align}
\begin{cases}
u \mbox{ is a maximal lifespan solution on } I;\\
u \mbox{ is spherically symmetric in space};\\
u \mbox{ satisfies the improved Duhamel formula \eqref{improve}}.
\end{cases}\label{condition}
\end{align}
By time translation invariance and without loss of generality we also assume $[0,1]\subset I$.

\subsection{Localization for kinetic energy}

The purpose of this subsection is to establish the uniform in time localization of the kinetic energy for solutions satisfying the condition \eqref{condition}. More precisely, we will prove
\begin{prop}[Kinetic energy localization] \label{miss}
Let $u$ satisfy \eqref{condition}. Then there exists a function $C(\eta)$ such that
\begin{align*}
\|\phi_{>C(\eta)} \nabla u(t)\|_{L_x^2}\le \eta, \ \ \forall \ \eta>0,\,\forall\ t \in I.
\end{align*}
\end{prop}

As shown in the proof of Theorem 1.14, Theorem 1.15 in \cite{lz:2d}
(see page 31), Proposition \ref{miss} will follow immediately from
the following two propositions which concern the decay of each
single frequency.

\begin{prop}[Frequency decay estimate]\label{fre_decay} Let $u$ satisfy \eqref{condition}.
 Let $\eps=\frac {d-1}d$. Then for any $t\in I$ and $N\ge 1$, we have
\begin{align}\label{comp}
\|\phi_{>1}P_N u(t)\|_{L_x^2}\lsm N^{-1-\eps}.
\end{align}
\end{prop}

\begin{rem}
The decay $N^{-1-\frac{d-1}{d}}$ may seem a bit surprising since the
exponent $1+\frac {d-1}d$ is bigger than the regularity of the
nonlinearity $1+\frac 4d$ for dimension $d>5$. However this is not
contradictory since in \eqref{comp} we are only considering the part
of the solution away from the origin. In this regime the additional
regularity of the solution comes from the smoothing effects of the
Schr\"odinger equation and the radial symmetry. On the other hand
for the part of the solution near the origin, we only obtain Sobolev
regularity $H^{s}$ for some $s<1+\frac 4d$ (see \eqref{an_bound}).
\end{rem}

\begin{prop}[Spatial decay estimate]\label{spa_decay} Let $u$ satisfy \eqref{condition}. Let $N_0, N_1$ be two dyadic numbers. Then there exist
 $R_0=R_0(N_0, N_1)$ and $\delta=\delta(d)$ such that for all $R\ge R_0$, $N\in [N_0, N_1]$ and $t\in I$, we have
\begin{align*}
\|\phi_{>R} P_N u(t)\|_{L_x^2}\lsm R^{-\delta}.
\end{align*}
\end{prop}

The proof of both propositions have been presented, in various
forms, in \cite{klvz} and \cite{lz:2d}. We sketch the proofs here
for the sake of completeness. The proof of Proposition \ref{miss}
will be skipped since it follows directly from Proposition
\ref{fre_decay} and Proposition \ref{spa_decay}.

\vspace{0.3cm}

{\it{Proof of Proposition \ref{fre_decay}:}}
\begin{proof}
We first use the in/out decomposition and triangle inequality to bound
\begin{align*}
\|\phi_{>1}P_N u(t)\|_2 \le \|\phi_{>1} P_N^+ u(t)\|_2+\|\phi_{>1}P_N^- u(t)\|_2.
\end{align*}
Since the two terms give the same contribution, we only estimate, for instance, the outgoing piece. For this piece, we use forward Duhamel formula. Moreover, we will split the integral into different time regimes and introduce the spatial cutoffs. We have
\begin{align}
\|\phi_{>1}P_N^+ u(t)\|_2&\lsm \|\phi_{>1} P_N^+\int_t^{\sup I} e^{i(t-s)\Delta} F(u(s)) ds\|_2\notag\\
&\lsm \|\phi_{>1}P_N^+\int_0^{\sup I-t} e^{-is\Delta}F(u(t+s)) d\tau\|_2\notag\\
&\lsm \|\phi_{>1} P_N^+\int_0^{\frac 1N} e^{-is\Delta}\phi_{>\frac 12}F(u(t+s))d s\|_2\label{main_s}\\
&\quad+ \|\phi_{>1}P_N^+\int_0^{\frac 1N} e^{-is\Delta}\phi_{\le \frac 12} F(u(t+s)) ds\|_2\label{tail_s}\\
&\quad+\|\phi_{>1}P_N^+ \int_{\frac 1N}^{\sup I-t} e^{-is\Delta}\phi_{>Ns/2}F(u(t+s)) ds\|_2\label{main_l}\\
&\quad+\|\phi_{>1}P_N^+ \int_{\frac 1N}^{\sup I-t} e^{-is\Delta}\phi_{\le Ns/2} F(u(t+s))d s\|_2.\label{tail_l}
\end{align}
The main contribution comes from \eqref{main_s} and \eqref{main_l}.
To estimate \eqref{main_s}, we drop the bounded operator
$\phi_{>1}P_N^+$ and commute the frequency cutoff $\tilde P_N$ with
the spatial cutoff $\phi_{>1}$(this produces a harmless high order
term by the mismatch estimate Lemma \ref{L:mismatch_fre}). Thus we
have
\begin{align}
\eqref{main_s}\lsm \|\int_0^{\frac 1N} e^{-is\Delta}\phi_{>\frac 12}P_{N/8<\cdot\le 8N}
F(\phi_{>\frac 14}u(t+s)) ds\|_2+N^{-10}.\label{c1}
\end{align}
%Now we need a crucial result from \cite{lz:2d}. Since $M(u)=M(Q)$ and $u\in H_x^1$, Lemma ?? in \cite{lz:2d} gives that
%\begin{align}\label{uniform}
%\|\phi_{>c}\nabla u(t)\|_2 \lsm_c 1.
%\end{align}
We now use weighted Strichartz Lemma \ref{L:wes} to estimate the
last term as
\begin{align*}
\eqref{main_s}&\lsm \|P_{N/8<\cdot \le 8N}F(\phi_{>\frac 14}u(t+s))\|_{L_s^{\frac d{d-1}}L_x^{\frac{2d}{d+4}}([0,\frac 1N])}+N^{-10}.\\
&\lsm N^{-\frac{d-1}d}.
\end{align*}
The estimate of \eqref{main_l} follows the similar way. Applying the
mismatch estimate and weighted  Strichartz inequality, we have
\begin{align*}
\eqref{main_l}&\lsm \|\int_{\frac 1N}^{\sup I-t} e^{-is\Delta}\phi_{>Ns/2}P_{N/8<\cdot\le 8N}
F(\phi_{>Ns/4}u(t+s)) ds\|_2+N^{-10}\\
&\lsm \|(Ns)^{-\frac{2(d-1)}d}P_{N/8<\cdot\le 8N}F(\phi_{>Ns/4}u(t+s))\|_{L_s^{\frac d{d-1}}L_x^{\frac{2d}{d+4}}([\frac 1N, \sup I-t))}+N^{-10}\\
&\lsm N^{-\frac{2(d-1)}d}\|s^{-\frac{2(d-1)} d}\| F(\phi_{>Ns/4}u(t+s))\|_{L_x^{\frac{2d}{d+4}}}\|_{L_s^{\frac d{d-1}}([\frac 1N, \sup I-t))}+N^{-10}\\
&\lsm N^{-\frac{d-1}d}.
\end{align*}
Finally we consider the contribution from the tail terms \eqref{tail_l} and \eqref{tail_s}. Applying Proposition \ref{P:P properties}, we bound the kernel as follows:
\begin{align*}
|(\phi_{> 1} P_N^+ e^{-is\Delta}\phi_{\le \frac 12})(x,y)|&\lsm N^{-9d}\langle N(x-y)\rangle^{-10 d}, \ \forall\ 0<s\le\frac 1N.\\
|(\phi_{>1}P_N^+ e^{-is\Delta}\phi_{\le Ns/2})(x,y)|&\lsm N^d \langle N^2 s\rangle^{-10 d}
\langle N(x-y)\rangle^{-10 d}\\
&\lsm N^{-9 d}\langle N(x-y)\rangle^{-10 d},\ \forall\ s>\frac 1N.
\end{align*}
The desired decay then follows from the kernel estimate and a simple use of Young's inequality.
Combining the estimates of these four pieces together, we obtain
\begin{align*}
\|\phi_{>1}P_N u(t)\|_{L_x^2}\lsm N^{-\frac {d-1}d}, \ \forall\ t\in I.
\end{align*}
Moreover it is easy to check that, after notational change, the same analysis establishes
\begin{align}\label{first}
\|\phi_{>c} P_N u(t)\|_{L_x^2}\lsm_c N^{-\frac{d-1}d},\ \forall\ t\in I.
\end{align}
This implies
\begin{align}\label{aw}
\||\nabla|^{\frac{d-1}d-}(\phi_{>c}u(t))\|_{L_x^2}\lsm_c 1, \ \forall\ t\in I.
\end{align}
Now we can upgrade the decay \eqref{first} by inserting \eqref{aw} when we repeat the same argument as above. For example, using Bernstein, \eqref{aw}, \eqref{main_s} can be re-estimated as follows:
\begin{align*}
\eqref{main_s}&\lsm \|P_{N/8<\cdot\le 8N} F(\phi_{>\frac 14} u(t+s))\|_{L_s^{\frac d{d-1}}L_x^{\frac{2d}{d+4}}([0,\frac 1N])}\\
&\lsm N^{-2\frac d{d-1}+}\||\nabla|^{\frac {d-1}d-}F(\phi_{>\frac 14} u(t+s))\|_{L_s^{\infty}L_x^{\frac{2d}{d+4}}([0,\frac 1N])}\\
&\lsm N^{-\frac{2(d-1)}d+}.
\end{align*}
The same computation applies to \eqref{main_l}, so we get
\begin{align*}
\|\phi_{>c}P_N u(t)\|_2\lsm_c N^{-\frac{2(d-1)}d+},\ \forall\ t\in I.
\end{align*}
Another repetition of the argument yields \eqref{comp} for $\eps=\frac {d-1}d$.

\end{proof}

The proof of Proposition \ref{spa_decay} has the same spirit as the proof of Proposition \ref{fre_decay}. So here we only briefly sketch the proof.

{\it Proof of Proposition \ref{spa_decay}:}
\begin{proof} Using in/out decomposition, it suffices to consider the piece
\begin{align*}
\|\phi_{>R}P_N^+ u(t)\|_2,
\end{align*}
for which we use forward Duhamel formula to express $u(t)$. This further reduces our consideration to the following integral
\begin{align*}
\|\phi_{>R}P_N^+ \int_0^{\sup I -t} e^{-is\Delta }F(u(t+s))ds\|_2.
\end{align*}
Now we spit the time integral into regimes where $0<s<\frac R{100
N}$, and $s>\frac R{100 N}$. For the small time regime, we insert
the spatial cutoff $\phi_{>R/2}$ and $\phi_{\le R/2}$. For the large
time regime, we insert the spatial cutoff $\phi_{>Ns/2}$ and
$\phi_{\le Ns/2}$. As indicated in the proof of Proposition
\ref{fre_decay}, the pieces with cutoff near the origin will give
arbitrary decay in $R$ by using the kernel estimate Proposition
\ref{P:P properties}. The pieces with cutoff away from the origin
can be dealt with by the weighted Strichartz estimate. The point
here is that since the frequencies are fixed in the dyadic interval
$[N_0,N_1]$, we can take $R$ sufficiently large to cancel any $N$
dependent quantity.
\end{proof}

%\begin{rem}
%One can see Proposition 3.1 in \cite{klvz} and Proposition 4.1--4.2 in \cite{lz:2d} for the proofs of both propositions. In particular, by using a refined analysis as in \cite{lz:2d}, it is not difficult to check that $\eps$ can be taken to be $\frac{d-1}d$. Also it is worthwhile pointing out that here since $u$ satisfies the condition \eqref{condition}, one can apply the improved Duhamel formula \emph{both backward and forward in time}, hence the contribution of the initial data does not appear on the right hand side of \eqref{comp}.
%\end{rem}
%
%Combining these two together, we have the following
%
%\begin{prop}[Kinetic energy localization, \cite{lz:2d}] \label{miss}
%Let $u$ satisfy \eqref{condition}. Then there exists a function $C(\eta)$ such that
%\begin{align*}
%\|\phi_{>C(\eta)} \nabla u(t)\|_{L_x^2}\le \eta, \ \ \forall \ \eta>0,\, t \in I.
%\end{align*}
%\end{prop}

\subsection{Local iteration to prove $ H_x^{1+}$ regularity.}
 In this part, we prove $u(0)=u_0\in H_x^{1+}$. This amounts to showing $\|P_{\ge N}u_0\|_{L_x^2}\lsm N^{-1-}$ for $N\ge 1$. Using Proposition \ref{fre_decay}, we first show the quantity
$\|P_{\ge N} u_0\|_{L_x^2}$ is determined by the dual Strichartz norm of the nonlinearity on the local time interval $[0,\frac 1{\sqrt N}]$.

\begin{lem}\label{reduce}
Let $u$ satisfy \eqref{condition}. Let $\eps=\frac {d-1}d$. Then for any $N\ge 1$,  we have
\begin{align}\label{a1}
\|P_{\ge N} u_0\|_{L_x^2}\le C(d,\|u_0\|_{L_x^2})( N^{-1-\eps}
+\|P_{\ge N}F(u)\|_{L_{t,x}^{\frac{2(d+2)}{d+4}}([0,\frac
1{\sqrt{N}}]\times\R^d)}).
\end{align}
\end{lem}

\begin{rem}
Here the choice of the time interval cutoff at $N^{-\frac 12}$ is
not special. Perhaps a more natural choice is $\frac 1 N$ since the
solution propagate at speed $N$ and one is localizing to spacial
scale $O(1)$. This latter choice would also work for our iteration
scheme.
\end{rem}

\begin{proof} Since by Proposition \ref{fre_decay},  $\|\phi_{>1} P_{\ge N} u_0\|_{L_x^2}\lsm N^{-1-\eps}$,
we only need to estimate the piece $\|\phi_{\le 1}P_{\ge N} u_0\|_{L_x^2}.$ In
the following, the implicit constants are allowed to depend on $d$
and $\|u_0\|_{L_x^2}$. By the improved Duhamel formula we get
\begin{align}
\|\phi_{\le 1} P_{\ge N} u_0\|_{L_x^2}&\le \|\phi_{\le 1}P_{\ge
N}\int_0^{\sup I} e^{-i\tau\Delta}F(u(\tau))d\tau\|_{L_x^2}\notag\\
&\le \|\phi_{\le 1}P_{\ge N}\int_0^{\frac 1{\sqrt N}}
e^{-i\tau\Delta} F(u(\tau))d\tau\|_{L_x^2}\label{a2}\\
&\quad+\|\phi_{\le 1} P_{\ge N}\int_{\frac 1{\sqrt N}}^{\sup I}
e^{-i\tau\Delta}\phi_{\le {N\tau}/{8}}
F(u(\tau))d\tau\|_{L_x^2}\label{a3}\\
&\quad+\|\phi_{\le 1} P_{\ge N}\int_{\frac 1{\sqrt N}}^{\sup I}
e^{-i\tau\Delta} \phi_{>N\tau/8}
F(u(\tau))d\tau\|_{L_x^2}.\label{a4}
\end{align}
For \eqref{a2}, we use Strichartz to bound it by
\begin{align*}
\|P_{\ge N}F(u)\|_{L_{t,x}^{\frac{2(d+2)}{d+4}}(\insd)}.
\end{align*}
For \eqref{a3}, using the kernel estimate with $m=10d$, we have
\begin{align*}
\eqref{a3}&\le \sum_{M\ge N}\|\phi_{\le 1} P_M\int_{\frac 1{\sqrt
N}}^{\sup I} e^{-i\tau\Delta} \phi_{\le N\tau/8}
F(u(\tau))d\tau\|_{L_x^2}\\
&\lsm \sum_{M\ge N} M^{d-20d}\int_{\sqn}^{\sup I} \tau^{-10d}\|\langle
M|\cdot|\rangle^{-10d} *F(u)\|_{L_x^2}d\tau\\
&\lsm \sum_{M\ge N} M^{d-20d}M^{\frac
12(10d-1)}\|F(u)\|_{L_{\tau}^{\infty}L_x^{\frac{2d}{d+4}}}\|\langle
M|\cdot|\rangle^{-10d}\|_{L_x^{\frac d{d-2}}}\\
&\lsm \sum_{M\ge N}{M^{\frac 32(1-10d)}}\\
&\lsm N^{-10}.
\end{align*}
For \eqref{a4}, by the triangle inequality, we have
\begin{align}
\eqref{a4}&\lsm \|P_{\ge N}\int_{\sqn}^{\sup I} e^{-i\tau\Delta}
\phi_{>N\tau/8}F(u\phi_{>1/8})(\tau)d\tau\|_{L_x^2}\notag\\
&\lsm \|P_{\ge N}\int_{\sqn}^{\sup I}
e^{-i\tau\Delta}\phi_{>N\tau/8}
P_{\le N/8}F(u\phi_{>1/8})(\tau)d\tau\|_{L_x^2}\label{a4a}\\
&\quad +\|P_{\ge N}\int_{\sqn}^{\sup I}
e^{-i\tau\Delta}\phi_{>N\tau/8} P_{>N/8}F(u\phi_{>1
/8})(\tau)d\tau\|_{L_x^2}\label{a4b}
\end{align}
For the term \eqref{a4a}, we use the mismatch estimate Lemma
\ref{L:mismatch_fre} and Bernstein to bound it as
\begin{align*}
\eqref{a4a}&\lsm \int_{\frac 1{\sqrt N}}^{\sup I} (N^2 \tau)^{-10d}\|P_{\le \frac N8} F(u\phi_{>\frac 18})\|_{L_x^2} d\tau\\
&\lsm \int_{\sqn}^{\sup I}(N^2\tau)^{-10d} N^2 d\tau\\
&\lsm N^{-5}.
\end{align*}
For the term \eqref{a4b}, we use weighted Strichartz to estimate and Proposition \ref{miss} to get
\begin{align*}
\eqref{a4b}&\lsm \|(N\tau)^{-\frac{2(d-1)}d}P_{>N/8}
F(u\phi_{>1/8})\|_{L_{\tau}^{\frac
d{d-1}}L_x^{\frac{2d}{d+4}}([\sqn,\sup I)\times\R^d)}\\
&\lsm N^{-\frac{2(d-1)}d}\|\tau^{-\frac{2(d-1)}d}\|_{L_{\tau}^{\frac
d{d-1}}([\sqn,\sup I))}\cdot N^{-1}\|\nabla P_{>N/8}
F(u\phi_{>1/8})\|_{L_{\tau}^{\infty}L_x^{\frac{2d}{d+4}}}\\
&\lsm N^{-1-\frac{3(d-1)}{2d}}.
\end{align*}
This finishes the proof of Lemma \ref{reduce}.

\end{proof}

Now we further estimate the dual Strichartz norm of the nonlinearity.

\begin{lem}[Dual Strichartz norm control]\label{L:dua_control} Let $u$ satisfy \eqref{condition}.
Let $\beta>0$, $N_0\ge 1$, $N>\frac 1{\beta}N_0$. Then for any
$0<s<1+\frac 4d$, we have
\begin{align}
&\|P_{\ge N} F(u)\|_{L_{t,x}^{\frac{2(d+2)}{d+4}}(\insd)} \notag \\
\lsm
& \|u\|_{S(\ins)}^{\frac 4d}\sum_{M\le \beta N}\biggl(\frac MN\biggr)^s\|P_M u\|_{S(\ins)}\notag\\
&\quad+\|u_{>\beta N}\|_{S(\ins)}\biggl(N_0^{\frac 4{d+2}}N^{-\frac 1{d+2}}
+\|u_{>N_0}\|_{L_{\tau}^\infty L_x^2}^{\frac
8{d(d+2)}}\|u_{>N_0}\|_{S(\ins)}^{\frac
4{d+2}}\biggr).\label{dual_norm}
\end{align}
\end{lem}
\begin{proof}
By splitting $u$ into low, medium and high frequencies:
\begin{align*}
u=u_{\le N_0}+u_{N_0<\cdot\le\beta N}+u_{>\beta N},
\end{align*}
we write
\begin{align}\label{fu}
F(u)=F(u_{\le \beta N})+O(u_{>\beta N}|u_{\le N_0}|^{\frac 4d})
+O(u_{>\beta N }|u_{>N_0}|^{\frac 4d}).
\end{align}
The contribution due to the first term can be estimated as follows.
By using Lemma \ref{special}, we have
\begin{align*}
&\qquad\|P_{\ge N} F(u_{\le \beta
N})\|_{L_{t,x}^{\frac{2(d+2)}{d+4}}(\insd)}\\
&\lsm N^{-s}
\||\nabla|^s
P_{\ge N} F(u_{ \le \beta N})\|_{L_{t,x}^{\frac{2(d+2)}{d+4}}(\insd)}\\
&\lsm N^{-s}\||\nabla|^s u_{\le \beta
N}\|_{L_{t,x}^{\frac{2(d+2)}d}(\insd)}\|u_{\le \beta
N}\|_{L_{t,x}^{\frac{2(d+2)}d}(\insd)}^{\frac 4d}\\
&\lsm \|u\|_{S(\ins)}^{\frac 4d}\sum_{M\le \beta N}\biggl(\frac
MN\biggr)^s \|P_M u\|_{S(\ins)}.
\end{align*}
For the contribution due to the second part of \eqref{fu}, we use
Bernstein to get
\begin{align*}
&\qquad\|u_{>\beta N} |u_{\le N_0}|^{\frac
4d}\|_{L_{t,x}^{\frac{2(d+2)}{d+4}}(\insd)}\\
&\lsm \|u_{>\beta N}\|_{L_{t,x}^{\frac{2(d+2)}d}(\insd)}\|u_{\le
N_0}\|_{L_{t,x}^{\frac{2(d+2)}d}(\insd)}^{\frac 4d}\\
&\lsm \|u_{>\beta N}\|_{S(\ins)}N_0^{\frac 4{d+2}} N^{-\frac
1{d+2}}\|u_{\le N_0}\|_{L_\tau^{\infty}L_x^2}^{\frac 4d}\\
&\lsm \|u_{>\beta N}\|_{S(\ins)}N_0^{\frac 4{d+2}}N^{-\frac 1{d+2}}.
\end{align*}
For the third term in \eqref{fu}, we use H\"older and interpolation to get
\begin{align*}
&\quad \|u_{>\beta N} |u_{> N_0}|^{\frac
4d}\|_{L_{t,x}^{\frac{2(d+2)}{d+4}}(\insd)}\\
&\lsm \|u_{>\beta N}\|_{L_{t,x}^{\frac{2(d+2)}d}(\insd)}\|u_{>
N_0}\|_{L_{t,x}^{\frac{2(d+2)}d}(\insd)}^{\frac 4d}\\
&\lsm \|u_{>\beta N}\|_{S(\ins)}\|u_{>
N_0}\|_{L_t^{\infty}L_x^{2}(\insd)}^{\frac 8{d(d+2)}}\|u_{>
N_0}\|_{S(\ins)}^{\frac 4{d+2}}.
\end{align*}
Collecting the three pieces together, we get \eqref{dual_norm}.
\end{proof}

Now by Strichartz estimate,
\begin{align*}
\|P_{\ge N} u\|_{S(\ins)}\lsm \|P_{\ge N} u_0\|_{L_x^2}+\|P_{\ge N} F(u)\|_{L_{t,x}^{\frac{2(d+2)}{d+4}}(\insd)}
\end{align*}
and the latter is in turn determined by $\|P_{\ge N} u\|_{S(\ins)}$ due to Lemma \ref{reduce} and Lemma \ref{L:dua_control}. This enables us to set up a recurrent relation for $\|P_{\ge N} u\|_{S(\ins)}$.

We define
\begin{align*}
A_N=\|P_{\ge N} u\|_{S(\ins)}.
\end{align*}

Since locally the Strichartz norm of $u$ is bounded, we can write
\begin{align*}
A:=\|u\|_{S([0,1])}+1<\infty.
\end{align*}
Using Strichartz inequality, Lemma \ref{reduce}, Lemma
\ref{L:dua_control} and taking $s=1+\frac 2d$, we obtain
\begin{align}
A_N&\le C(d)(\|P_{\ge N}u_0\|_{L_x^2}+\|P_{\ge N}
F(u)\|_{L_{t,x}^{\frac{2(d+2)}{d+4}}(\insd)})\notag\\
& \le C(d,\|u_0\|_{L_x^2})\biggl( N^{-1-\eps}+\notag\\
&\qquad A^{\frac 4d}\sum_{M\le \beta N}\biggl(\frac
MN\biggr)^{1+\frac 2d}
\|P_M u\|_{S(\ins)}\label{b1}\\
&\qquad+\|P_{\ge \beta N}u\|_{S(\ins)}(N_0^{\frac 4{d+2}}N^{-\frac
1{d+2}}+ A^{\frac 4{d+2}}\|u_{\ge
N_0}\|_{L_t^{\infty}L_x^2(\ins)}^{\frac
8{d(d+2)}})\biggr).\label{b2}
\end{align}
For \eqref{b1}, we do a little modification. Noting $P_M=P_MP_{\ge
M/2}$, we have
\begin{align*}
\eqref{b1}&\lsm A^{\frac 4d}\sum_{M\le \beta N}\biggl(\frac
MN\biggr)^{1+\frac 2d}\|P_{\ge M/2}u\|_{S(\ins)}\\
&\lsm A^{\frac 4d}\sum_{M\le 2\beta N}\biggl( \frac
MN\biggr)^{1+\frac 2d}\|P_{\ge M}u\|_{S(\ins)}.
\end{align*}
We shall take $\beta$ to be sufficiently small. The constraint on
$\beta$ will be specified later.

Now we absorb \eqref{b2} into \eqref{b1} through taking suitable
parameters. First we take $N_0=N_0(\beta,A)$ such that
\begin{align*}
A^{\frac 4{d+2}}\|u_{>N_0}\|_{L_t^{\infty}L_x^2([0,1])}^{\frac
8{d(d+2)}}\le \frac 1{100} \beta^{1+\frac 2d}.
\end{align*}
This is certainly possible since $u\in C([0,1],\, L_x^2)$ and $[0,1]$ is a compact interval. Then we assume $N\ge
M_0$ where
\begin{align}\label{assum_m0}
M_0^{-\frac 1{d+2}}N_0^{\frac 4{d+2}}\le \frac 1{100}\beta^{1+\frac
2d}.
\end{align}
Under these restrictions we have
\begin{align}
\eqref{b2}\le \frac 12\beta^{1+\frac 2d}\|P_{\ge\beta N}
u\|_{S(\ins)}.
\end{align}
Therefore we get for all $N\ge M_0$ that
\begin{align*}
A_N&\le C(d,\|u_0\|_{L_x^2}) \biggl( N^{-1-\eps}+ \sum_{M\le 2\beta
N}\biggl(\frac MN\biggr)^{1+\frac 2d}\|P_{\ge M}
u\|_{S(\ins)}\biggr)\\
&\le N^{-1-\frac \eps 2}+\sum_{M\le 2\beta N}\biggl(\frac
MN\biggr)^{1+\frac 1d}\|P_{\ge M}u\|_{S(\ins)},
\end{align*}
where in the last inequality we have killed the constant
$C(d,\|u_0\|_{L_x^2})$. This is possible by first taking $\beta$
sufficiently small, then taking $M_0$ large enough.

Now, we split the summation into $M\le M_0$ and $M>M_0$. For large
$M$, we trivially bound the summand by
\begin{align*}
\biggl(\frac MN\biggr)^{1+\frac 1d}A_M.
\end{align*}
Then we sum all the pieces for small $M$, this  gives that
\begin{align*}
\sum_{M\le M_0}\biggl(\frac MN\biggr)^{1+\frac 1d}\|P_{\ge
M}u\|_{S(\ins)}\lsm AM_0^{1+\frac 1d} N^{-1-\frac 1d}.
\end{align*}
Finally we establish the following recurrence relation for $A_N$:
Let $s=\frac 1d+1$. Then there exists $C_1>0$ such that for all
$N\ge M_0$,
\begin{align}\label{rec_rel}
A_N\le C_1 M_0^s N^{-s}+\sum_{M_0< M\le 2\beta N}\biggl(\frac
MN\biggr)^s A_M
\end{align}
This combined with the trivial bound $A_N\le A$ will give us the
final control on $A_N$
\begin{align}
A_N\le C(A, M_0) N^{-s+}, \ \forall\ N\ge M_0,\label{an_bound}
\end{align}
if we apply the following lemma:

\begin{lem}[Recursive control]\label{general}
Let $s>1$, $\gamma>0$ and $s-\gamma>1$. Let $C_1>0$ be such that
for all
$N\ge M_0$,
\begin{align}
A_N&\le C_1 M_0^s N^{-s}+\sum_{M_0\le M\le \beta^\prime N}\biggl(\frac
MN\biggr)^s
A_M,\label{an_control}\\
A_N&\le A.\label{trivial}
\end{align}
Then there exists a constant $c(s,\gamma,A)>0$ such that for
all $0<\beta^\prime<c(s,\gamma,A)$, we have
\begin{align}
A_N\le 2C_1 M_0^s N^{-s+\gamma},\ \forall\, N\ge M_0.\label{final_an}
\end{align}
\end{lem}
\begin{proof}
We will inductively prove
\begin{align}\label{induct}
A_N\le 2C_1 M_0^s N^{-s+\gamma}+(\beta^\prime)^j.
\end{align}
First, plugging the bound \eqref{trivial} into \eqref{an_control},
we get
\begin{align*}
A_N\le C_1M_0 N^{-s}+C(s)A(\beta^\prime)^s\le 2C_1 M_0 N^{-s+\gamma} +\beta^\prime.
\end{align*}
by requiring $(\beta^\prime)^{s-1}<\frac 1{100C(s)A}$. This establishes
\eqref{induct} for $j=1$.

Now assume \eqref{induct} hold for $j$-th step, we plug this
bound into \eqref{an_control} to compute
\begin{align*}
A_N&\le C_1 M_0^s N^{-s}+2C(s)(\beta^\prime)^{\gamma}\cdot C_1M_0^s
N^{-s+\gamma}+C(s)(\beta^\prime)^{s-1}\cdot(\beta^\prime)^{j+1}\\
&\le 2C_1M_0^s N^{-s+\gamma}+(\beta^\prime)^{j+1}.
\end{align*}
by requiring $(\beta^\prime)^{\gamma}<\frac 1{100 C(s)}$. This establishes
\eqref{induct} for $j+1$.

Finally, \eqref{final_an} follows by taking $j\to\infty$ in
\eqref{induct}.
\end{proof}

\end{document}